\documentclass[leqno]{article}
\usepackage{amsmath,amssymb,amsthm}

\usepackage{titlesec,hyperref}

\usepackage{color}

\usepackage{fancyhdr}
\titleformat{\subsection}{\it}{\thesubsection.\enspace}{1pt}{}

\newtheorem{theo}{Theorem}[section]
\newtheorem{lemm}[theo]{Lemma}

\newtheorem{prop}[theo]{Proposition}
\newtheorem{rema}[theo]{Remark}
\numberwithin{equation}{section}

\allowdisplaybreaks

\begin{document}
\title{Global Well-posedness for the Generalized Navier-Stokes System
\hspace{-4mm}
}

\author{Zeng Zhang$^1$
\quad Zhaoyang Yin$^2$ \\[10pt]
Department of Mathematics, Sun Yat-sen University,\\
510275, Guangzhou, P. R. China.\\[5pt]
}
\footnotetext[1]{Corresponding author. Email: \it zhangzeng534534@163.com; Tel.:+8613725201644;\\
Address: No. 135, Xingang Xi Road, Guangzhou, 510275, P. R. China}
\footnotetext[2]{Email: \it mcsyzy@mail.sysu.com.cn}
\date{}
\maketitle

\begin{abstract}
In this paper we investigate well-posedness of the Cauchy problem
of the three dimensional generalized Navier-Stokes system. We
first establish local well-posedness of the GNS system for any
initial data in the Fourier-Herz space $\chi^{-1}$. Then we show
that if the $\chi^{-1}$ norm of the initial data is smaller than
C$\nu$ in the GNS system where $\nu$ is the viscosity coefficient,
the corresponding solution exists globally in time. Moreover, we
prove global well-posedness of the Navier-Stokes system without
norm restrictions on the corresponding solutions provided the
$\chi^{-1}$ norm of the initial data is less than $\nu.$ Our
obtained results cover and improve recent results in \cite{Zhen
Lei,wu}.

\vspace*{5pt}

\end{abstract}

\vspace*{10pt}


\tableofcontents

\section{Introduction}
In this paper we consider the generalized Navier-Stokes system
\begin{align}\tag{$GNS_\nu$}
\left\{
\begin{array}{ll}
\partial_tu-\nu \triangle u =Q(u,u),\quad t>0,\,x \in \mathbb{R}^3,\\[1ex]
u(x,0)=u_0,
\end{array}
\right.
\end{align}
with $\nu$ being a positive constant and the bilinear operator $Q$
defined as
\begin{align}\label{(1.1)}
Q^j(u,v)=\sum_{k,l,m=1}^3q^{j,m}_{k,l}\partial_m(u^kv^l),\,\,j=1,2,3,
\end{align}
where $q^{j,m}_{k,l}(a)=\sum_{n,p=1}^3a^{j,m,p,n}_{k,l}{\mathcal{F}}^{-1}(\frac{\xi_n \xi_p}{|\xi|^2}\widehat{a}(\xi))$,
 and $a^{j,m,p,n}_{k,l}$ are real numbers.

It is obvious that the incompressible Navier-Stokes system
\begin{align}\tag{$NS_{\nu}$}
\left\{
\begin{array}{ll}
\partial_tu-\nu \triangle u =\mathcal{P}\nabla(u\otimes u), \quad t>0,\,x \in \mathbb{R}^3,\\[1ex]
div~u=0,\\[1ex]
u(x,0)=u_0,
\end{array}
\right.
\end{align}
is a particular case of the system $(GNS_\nu)$. Here $u$ stands for the velocity field of the fluid, $\nu$ is the viscosity and $\mathcal{P}$ is the Leray projection operator defined by the formula:
\begin{align}
  \mathcal{F}(\mathcal{P}f)^j(\xi)=\mathcal{F}(f)^j(\xi)-\frac{1}{|\xi|^2}\sum_{k=1}^3\xi_j \xi_k\mathcal{F}(f)^k(\xi),\,j=1,2,3.
\end{align}
From now on  we always assume that the initial data $u_0$~is
divergence free and C denotes a generic constant.

It is well known that the space $BMO^{-1}$ is the largest space
which is included in the tempered distribution and enjoys the
property of translation and scaling invariant (see
\cite{Chemin,Koch} for instance). The global well-posedness for
the Navier-Stokes system in the space $BMO^{-1}$ was studied by
Koch and Tataru \cite{Koch}.  Many works in subspaces of the space
$BMO^{-1}$, such as $\dot{H}^{\frac{3}{2}-1}(\mathbb{R}^3)$,
$L^3(\mathbb{R}^3)$, and
$\dot{B}_{p,\infty}^{-1+\frac{3}{p}}(\mathbb{R}^3)$ also have been
down: the Navier-Stokes system is known to be globally well-posed
for sufficiently small date $u_0\in
\dot{H}^{\frac{3}{2}-1}(\mathbb{R}^3),$ and locally well-posed for
any $u_0\in \dot{H}^{\frac{3}{2}-1}(\mathbb{R}^3),$ as proved by
Fujita and Kato \cite{Fujita and Kato}; the global well-posedness
for small date is due to Kato \cite{Kato1} in $L^3(\mathbb{R}^3)$,
Cannone, Meyer and Planchon \cite{M. Cannone} in
$\dot{B}_{p,\infty}^{-1+\frac{3}{p}}(\mathbb{R}^3)$, respectively.
For more studies in these spaces, the reader may refer to
$\cite{2,1,4}.$ Let us mention that the well-posedness for the
Navier-Stokes system in those above spaces can be extended to the
generalized Navier-Stokes system \cite{Bahouri}.

 Recently Lei and Lin \cite{Zhen Lei} studied a new space $$\chi^{-1}=\{f\in \mathcal{D}'({\mathbb{R}^3});\int_{\mathbb{R}^3}|\xi|^{-1}|\widehat{f}|d\xi< \infty\},$$ which is contained in $BMO^{-1}$ \cite{Zhen Lei} and is equivalent to the Fourier-Herz space $\dot{\mathcal{B}}^{-1}_1$ \cite{wu}. They pointed out that $H^s(\mathbb{R}^3)\,(s>\frac{1}{2})\,\subseteq \chi^{-1},$ and they also presented an example to show that $H^{\frac{3}{2}-1}(\mathbb{R}^3)\nsubseteq \chi^{-1}.$ Here we give an example to show that $\chi^{-1}\nsubseteq \dot{H}^{\frac{3}{2}-1}(\mathbb{R}^3).$ In fact, let $f(x)= \mathcal{F}^{-1}(\frac{1}{|\xi|}h(|\xi|)),$ where $h(r)$ is defined by
 \begin{align*}
h(r)=\left\{
\begin{array}{ll}
  2^{\frac{j+1}{2}},&1-2^{-j}\leq r< 1-2^{-(j+1)},j=0,1,\cdots,\\[1ex]
  0,&r\geq1.
  \end{array}
\right.
 \end{align*}
 It is easy to deduce that
 \begin{align*}
   \int_{\mathbb{R}^3}|\xi|^{-1}||\widehat{f}(\xi)|d\xi&\leq C\int_0^\infty h(r)dr
   =C\sum_{j=0}^\infty 2^{\frac{j+1}{2}}2^{-(j+1)}< \infty.
 \end{align*}
 Similarly, we can also get
 \begin{align*}
   \int_{\mathbb{R}^3}|\xi||\widehat{f}(\xi)|^2d\xi&= C\int_0^\infty rh(r)^2dr
  \geq C\sum_{j=0}^\infty  (1-2^{-j})2^{j+1}2^{-(j+1)}= \infty.
 \end{align*}
 Combing the above two inequalities yields that $f\in \chi^{-1},$ but $f \not\in \dot{H}^{\frac{3}{2}-1}(\mathbb{R}^3).$ Thus, we conclude that $\chi^{-1}$ and $\dot{H}^{\frac{3}{2}-1}(\mathbb{R}^3)$ do not contain each other.

 Lei and Lin \cite{Zhen Lei} proved that if the initial data $u_0$ in $\chi^{-1}$ satisfying
 $\|u_0\|_{\chi^{-1}}<\nu,$ the system $(NS_\nu)$ admits a global mild solution.
 They also proved that this global mild solution is unique under the condition
 $\|u\|_{L^\infty(\mathbb{R}^{+};\chi^{-1})}< \nu.$ Cannone and Wu \cite{wu} gave a
 global well-posedness result for small initial data in a family of critical Fourier-Herz spaces
 $\dot{\mathcal{B}}^{-1}_q\,(q\in[1,2]).$ They also showed this global solution is unique under the
 condition $\|u\|_{L^\infty(\mathbb{R}^{+};\mathcal{B}^{-1}_q)\cap L^1(\mathbb{R}^{+};\mathcal{B}^{1}_q)}\leq \frac{50}{9}\|u_0\|_{\mathcal{B}^{-1}_q}$.  However, it is not clear whether there exists a solution to the system $(NS_\nu)$ for large initial data.
 Moreover, without norm restrictions on the solutions, is the uniqueness of solutions to the system $(NS_\nu)$ still valid?

In the paper we will give definite answers to these two questions.
In Section 3, we will solve  the system $(GNS_\nu)$ by means of a
contraction mapping argument ( see Lemma \ref{Main Lemma} below).
Thus, we can obtain a unique local mild solution to the system
$(GNS_\nu)$ for any initial data in $\chi^{-1}$ and prove the
corresponding solution will be global if the initial data is
sufficient small. Especially for the Navier-Stokes system, we show
that if $\|u_0\|_{\chi^{-1}}<\nu$ then the solution to the system
$(NS_\nu)$ will be unique and global without norm restrictions on
the solutions.
\begin{lemm}(\cite{Bahouri})\label{Main Lemma}  Let E be a Banach space, $\mathcal{B}$ a continuous bilinear map from $E \times E\rightarrow E,$ and a positive real number such that
$\alpha < \frac{1}{4\| \mathcal{B} \|},$ with
$$\| \mathcal{B}\|=\sup_{\|u\| \leq 1,\|v\| \leq 1}\|  \mathcal{B}(u,v)\|.$$
For any $a$ in the ball $B(0,\alpha)$ in $E$, then there exists a
unique $x$ in $B(0, 2\alpha)$ such that
$$x=a+\mathcal{B}(x,x).$$
\end{lemm}
We will also use the spaces
\begin{align}\chi^{i}=\{f\in \mathcal{S}'({\mathbb{R}^3});\int_{\mathbb{R}^3}|\xi|^{i}|\hat{f}|d\xi< \infty\},~~~i=-1, 0,1.\end{align}
The norm of $\chi^i$ is denoted by $\|\cdot\|_{\chi^i}.$ Let $T\in
(0,\infty].$ Note that the spaces $L^2([0,T];\chi^0)$ and
$L^\infty([0,T];\chi^{-1})\cap L^1([0,T];\chi^1)$ are Banach
spaces and are translation and shift invariant.

Now we are in the position to state our main results:
\begin{theo}\label{main theo}
Let $u_0$ be in $\chi^{-1}.$  There exists a positive time $T$
such that the system $(GNS_\nu)$ has a unique solution $u$ in
$L^2([0,T];\chi^0)$ which also belongs to
$$C([0,T];\chi^{-1})\cap L^1([0,T];\chi^1)\cap L^\infty([0,T];\chi^{-1}).$$
Let $T_{u_0}$ denote the maximal time of existence of such a solution. Then:\\
(i) There exists a constant $C$ such that if
$\|u_0\|_{\chi^{-1}}\leq C \nu ,$ then
$$T_{u_0}=\infty.$$
(ii) If $T_{u_0}$ is finite, then
$$\int^{T_{u_0}}_0 \|u(t)\|_{\chi^0}^2dt=\infty.$$

\end{theo}
\begin{theo}\label{main coro}
  Let $u_0$ be in $\chi^{-1}.$ There exists a positive time $T$ such that the system $(NS_\nu)$ has a unique solution u in $L^2([0,T];\chi^0)$ which also belongs to
$$C([0,T];\chi^{-1})\cap L^1([0,T];\chi^1)\cap L^\infty([0,T];\chi^{-1}).$$
Let $T_{u_0}$ denote the maximal time of existence of such a solution. Then:\\
(i) There exists a constant $C$ such that if $\|u_0\|_{\chi^{-1}}<
\nu,$ then
$$T_{u_0}=\infty.$$
(ii) If $T_{u_0}$ is finite, then
$$\int^{T_{u_0}}_0 \|u(t)\|_{\chi^0}^2dt=\infty.$$
\end{theo}

\begin{rema}
Although the result (i) in Theorem \ref{main theo} is the same as
Theorem (1.1) in \cite{Zhen Lei}, our method here is different
from their method in \cite{Zhen Lei}. In particular, our proof
relies on the obtained local well-posedness result and the blow-up
criterion (ii), but not on additional norm restrictions on the
corresponding solutions.
\end{rema}
\section{Preliminaries}
Let $B(u,v)$ be the solution to the heat equation
\begin{align}
\left\{
\begin{array}{l}
\partial_t B(u,v)-\nu \triangle B(u,v)=Q(u,v),  \\[1ex]
B(u,v)|_{t=0}=0,
\end{array}
\right.
\end{align}
with the bilinear operator $Q$ defined as in (\ref{(1.1)}).\\
Solving $(GNS_\nu)$ amounts to finding a fixed point for the map
$$u\mapsto e^{t\nu \triangle} u_0+B(u,u).$$
By Duhamel's formula in Fourier space and (\ref{(1.1)}), we have
\begin{align}
|\widehat{B(u,v)}(t,\xi)| &=
|\int_0^t e^{-\nu(t-s)|\xi|^2}\widehat{Q(u,v)}(s,\xi)ds|\\
\nonumber&\leq C\int_0^t e^{-\nu(t-s)|\xi|^2}|\xi| (|\widehat{u}|\ast |\widehat{v}|)(s,\xi)ds.
\end{align}
We now give two useful propositions which will be used in the sequel.\\
\begin{prop}\label{chazhi}
\begin{align}\|u\|_{L^2([0,T];\chi^{0})}^2\leq\|u\|_{L^\infty([0,T];\chi^{-1})}\|u\|_{L^1([0,T];\chi^{1})}.\end{align}
\end{prop}
{\noindent\bf Proof.} It is easy to check that
\begin{align*}
\|u\|_{L^2([0,T];\chi^{0})}^2&=\int^T_0 \left( \int_{\mathbb{R}^3} |\widehat{u}|(s,\xi) d\xi \right) ^2 dt\\
\nonumber&\leq \int^T_0 \left( \int_{\mathbb{R}^3}|\xi|^{-1} |\widehat{u}|(s,\xi)d\xi \right)
                    \left( \int_{\mathbb{R}^3}|\xi||\widehat{u}|(s,\xi) d\xi \right) dt\\
\nonumber&\leq\|u\|_{L^{\infty}([0,T];\chi^{-1})}\|u\|_{L^1([0,T];\chi^{1})}.
\end{align*}
\begin{prop}\label{B}
A constant C exists such that
\begin{align}
 \|B(u,v)\|_{L^2([0,T];\chi^{0})}\leq \frac{C}{\nu ^{\frac{1}{2}}}\|u\|_{L^2([0,T];\chi^{0})}\|v\|_{L^2([0,T];\chi^{0})}.
\end{align}
\end{prop}
{\noindent\bf Proof.} Thanks to the inequality $(2.2)$ and Minkowski's inequality, we have
\begin{align*}
\|B(u,v)\|_{L^2([0,T];\chi^{0})} &\leq C\left\|\int_{R^3} \int^T_0 I_{[0,t]}(s)e^{-\nu (t-s)|\xi|^2}|\xi| (|\widehat{u}|\ast |\widehat{v}|)(s,d\xi)ds d\xi \right\|_{L^2(0,T)}\\
\nonumber &\leq C\int_{R^3} \int^T_0 \left(\int^T_0 I_{[0,t]}(s)e^{-2\nu(t-s)|\xi|^2}|\xi|^2dt \right)^{\frac{1}{2}} (|\widehat{u}|\ast |\widehat{v}|)(s,d\xi)ds d\xi\\
\nonumber &\leq C\int_{R^3} \int^T_0 \frac{1}{\nu ^{\frac{1}{2}}}(|\widehat{u}|\ast |\widehat{v}|)(s,d\xi)ds d\xi\\
\nonumber &\leq \frac{C}{\nu ^{\frac{1}{2}}}\int^T_0\|\widehat{u(s)}\|_{L^1}\|\widehat{v(s)}\|_{L^1}ds\\
\nonumber &\leq \frac{C}{\nu ^{\frac{1}{2}}}\|u\|_{L^2([0,T];\chi^{0})}\|v\|_{L^2([0,T];\chi^{0})}.
\end{align*}

\section{Proofs of main theorems}
To prove the first part of Theorem \ref{main theo}. we shall use Lemma 1.1. Given some $u_0 \in \chi^{-1},$ thanks to Minkowski's inequality, we have
\begin{align}
\|e^{\nu t\triangle }u_0\|_{L^2([0,T];\chi^0)}&=\left(\int_0^T(\int_{R^3}e^{-\nu t|\xi|^2}|\widehat{u_0}|(\xi)d\xi)^2 dt\right)^{\frac{1}{2}}\\
\nonumber&\leq \int_{R^3}\left(\int_0^T e^{-2\nu t|\xi|^2}|\widehat{u_0}|^2(\xi)ds\right)^{\frac{1}{2}}d\xi\\
\nonumber&\leq \int_{R^3}\frac{1}{(2\nu |\xi|^2)^{\frac{1}{2}}}|\widehat{u_0}|(\xi)d\xi\\
 \nonumber&\leq \frac{1}{(2\nu)^{\frac{1}{2}}} \|u_0\|_{\chi^{-1}}.
\end{align}
Thus, combining Proposition \ref{B} and the inequality (3.1) gives that if $\|u_0\|_{\chi^{-1}}\leq \frac{\nu}{2^{\frac{3}{2}}C_0},$ with $C_0>C,$ then
$$\|e^{\nu t \triangle  }u_0\|_{L^2([0,T];\chi^0)}\leq \frac{1}{4\frac{C_0}{\nu^{\frac{1}{2}}}}<\frac{1}{4\|B\|}.$$
According to Lemma 1.1, there exists a unique solution of the system $(GNS_\nu)$ in the ball with center $0$
and radius $\frac{\nu^{\frac{1}{2}}}{2C_0}$ in the space $L^2([0,T];\chi^0).$\\
\vspace*{1ex}

We now consider the case of a large initial date $u_0\in \chi^{-1}.$ We shall split $u_0$ into a small part in $\chi^{-1}$ and a large part with compactly supported Fourier transform. For that, we fix some positive real number $\rho_{u_0}$ such that
\begin{align}
  \int_{|\xi|\geq \rho_{u_0}}|\xi|^{-1}|\widehat{u_0}|(\xi)d\xi\leq \frac{\nu}{2^{\frac{5}{2}}C_0}.
\end{align}
Using the inequality (3.1) again and defining $u_0^\flat=\mathcal{F}^{-1}(I_{B(0,\rho_{u_0})}(\xi)\widehat{u_0}(\xi)),$ we get
\begin{align*}
\|e^{\nu t \triangle  }u_0\|_{L^2([0,T];\chi^0)}\leq \frac{\nu^{\frac{1}{2}}}{8C_0}
 +\|e^{\nu t \triangle  }u_0^\flat\|_{L^2([0,T];\chi^0)}.
 \end{align*}
From which we can deduce that
 \begin{align*}
\|e^{\nu t\triangle }u_0^\flat\|_{L^2([0,T];\chi^{-1})}
&=\left(\int^T_0 (\int_{|\xi|\leq \rho_{u_0}}e^{-\nu t|\xi|^2} |\widehat{u_0}|(\xi) d\xi)^2dt\right)^{\frac{1}{2}}\\
\nonumber&=\left(\int^T_0 (\int_{|\xi|\leq \rho_{u_0}}e^{-\nu t|\xi|^2}|\xi||\xi|^{-1} |\widehat{u_0}|(\xi) d\xi)^2dt\right)^{\frac{1}{2}}\\
\nonumber& \leq \rho_{u_0}T^{\frac{1}{2}}\|u_0\|_{\chi^{-1}}.
\end{align*}
Thus if
\begin{align}
  T\leq (\frac{\nu^{\frac{1}{2}}}{8\rho_{u_0}C_0\|u_0\|_{\chi^{-1}}})^2,
\end{align}
then we conclude the existence of a unique solution in the ball with center 0 and radius $\frac{\nu^{\frac{1}{2}}}{2C_0}$ in the space $L^2([0,T];\chi^0).$\\

Next we claim that if $u$ is a solution of the system $(GNS_\nu)$ in $L^2([0,T];\chi^0),$ then $u$ also belongs to  $C([0,T];\chi^{-1})\cap L^1([0,T];\chi^1)\cap L^\infty([0,T];\chi^{-1}).$ In fact, it is easy to deduce that
\begin{align*}
\|B(u,v)(t)\|_{\chi^{-1}}
&\leq C \int_{\mathbb{R}^3}\int_0^t |\xi|^{-1}e^{-\nu (t-s)|\xi|^2}|\xi|(|\hat{u}|\ast |\hat{v}|)(s,\xi)dsd\xi\\
\nonumber&\leq C \int_{\mathbb{R}^3}\int_0^t e^{-\nu (t-s)|\xi|^2}(|\hat{u}|\ast |\hat{v}|)(s,\xi)dsd\xi\\
\nonumber&\leq C \int_0^t\int_{\mathbb{R}^3}(|\hat{u}|\ast |\hat{v}|)(s,\xi)d\xi ds\\
\nonumber&\leq C\|u\|_{L^2([0,T];\chi^0)}\|v\|_{L^2([0,T];\chi^0)}.
\end{align*}
Similarly, we have
\begin{align*}
\|B(u,v)\|_{L^1([0,T];\chi^{1})} &\leq C\int_0^T\int_{\mathbb{R}^3} \int^T_0 |\xi|I_{[0,t]}(s)e^{-\nu (t-s)|\xi|^2}|\xi| (|\hat{u}|\ast |\hat{v}|)(s,\xi)ds d\xi dt\\
\nonumber &\leq C\int_{\mathbb{R}^3} \int^T_0 \left(\int^T_0 I_{[0,t]}(s)e^{-\nu(t-s)|\xi|^2}|\xi|^2dt \right) (|\hat{u}|\ast |\hat{v}|)(s,\xi)ds d\xi\\
\nonumber &\leq C\int_{\mathbb{R}^3} \int^T_0 \frac{1}{\nu }(|\hat{u}|\ast |\hat{v}|)(s,\xi)ds d\xi\\
\nonumber &\leq \frac{C}{\nu }\|u\|_{L^2([0,T];\chi^0)}\|v\|_{L^2([0,T];\chi^0)}.
\end{align*}
Combing the above two inequalities yields that
$$B(u,u)\in L^\infty([0,T];\chi^{-1})\cap L^1([0,T];\chi^1).$$
Noticing  the following two facts:
\begin{align*}
\|e^{\nu t\triangle }u_0\|_{\chi^{-1}}\leq \int_{\mathbb{R}^3}e^{-\nu t|\xi|^2}|\xi|^{-1} |\widehat{u_0}|(\xi) d\xi \leq \|u_0\|_{\chi^{-1}},
\end{align*}
\begin{align*}
\|e^{\nu t\triangle }u_0\|_{L^1([0,T];\chi^{-1})}\leq  \int_0^T\int_{\mathbb{R}^3}e^{-\nu t|\xi|^2}|\xi| |\widehat{u_0}|(\xi) d\xi dt \leq \frac{1}{\nu} \|u_0\|_{\chi^{-1}},
\end{align*}
 we have
$$e^{\nu t\triangle }u_0 \in L^\infty([0,T];\chi^{-1})\cap L^1([0,T];\chi^1).$$
We can thus conclude that
$$u \in  L^\infty([0,T];\chi^{-1})\cap L^1([0,T];\chi^1).$$

To get further the regularity of $u(t,x)$ with respect to $t,$ we come back to the system $(GNS_\nu).$ It is obvious that $\triangle u$ is in $L^1([0,T];\chi^{-1}).$ Thanks to Proposition \ref{chazhi} we have
\begin{align}
  &\|Q(u,v)\|_{L^1([0,T];\chi^{-1})}\\
\nonumber=&\int_0^T \int_{\mathbb{R}^3}|\xi|^{-1}|\widehat{Q(u,v)}|(s,\xi)d\xi ds\\
\nonumber\leq &C\int_0^T \int_{\mathbb{R}^3}  (|\hat{u}|\ast |\hat{v}|)(s,\xi)d\xi ds\\
\nonumber\leq &C\|u\|_{L^2([0,T];\chi^{0})}\|v\|_{L^2([0,T];\chi^{0})}\\
\nonumber\leq &C\|u\|_{L^{\infty}([0,T];\chi^{-1})}^{\frac{1}{2}}\|u\|_{L^1([0,T];\chi^{1})}^{\frac{1}{2}}
\|v\|_{L^{\infty}([0,T];\chi^{-1})}^{\frac{1}{2}}\|v\|_{L^1([0,T];\chi^{1})}^{\frac{1}{2}}.
\end{align}
Thus, $Q(u,u) \in L^1([0,T];\chi^{-1}).$ Then we conclude that $\partial_t u \in L^1([0,T];\chi^{-1}).$ We hence prove the announced properties:
$$u \in C([0,T];\chi^{-1})\cap L^1([0,T];\chi^1)\cap L^\infty([0,T];\chi^{-1}).$$
\vspace*{1ex}

The uniqueness part relies on the following lemma.
\begin{lemm}\label{zhongyao}
Let $v$ be a solution in $C([0,T];\mathcal{S}'(R^3))$ of the Cauchy problem
\begin{align}
\left\{
\begin{array}{l}
\partial_t v-\nu \triangle v=f,  \\[1ex]
v|_{t=0}=v_0,
\end{array}
\right.
\end{align}
with $f\in L^1([0,T];\chi^{-1})$ and $v_0\in \chi^{-1}.$ Then for any $0\leq t_0\leq T$, we have
\begin{align}
\int_{\mathbb{R}^3} \sup_{0\leq t\leq t_0}|\hat{v}|(t,\xi)|\xi|^{-1}d\xi
\leq \|v_0\|_{\chi^{-1}}+\|f\|_{L^1([0,t_0];,\chi^{-1})},
\end{align}
\begin{align}
  \|v\|_{L^\infty([0,t_0];\chi^{-1})}+\nu \|v\|_{L^1([0,t_0];\chi^{1})}
  \leq 2(\|v_0\|_{\chi^{-1}}+\|f\|_{L^1([0,t_0];\chi^{-1})}).
\end{align}
\end{lemm}
{\noindent\bf Proof.}
Due to Duhamel's formula in Fourier space, we can write that
\begin{align*}
\widehat{v}(t,\xi) =e^{-\nu t|\xi|^2}\widehat{v_0}(\xi)+\int_0^t e^{-\nu(t-s)|\xi|^2}\widehat{f}(s,\xi)ds.
\end{align*}
Then  we have
\begin{align*}
  |\widehat{v}(t,\xi)| \leq e^{-\nu t|\xi|^2}|\widehat{v_0}|(\xi)+\int_0^t e^{-\nu(t-s)|\xi|^2}|\widehat{f}|(s,\xi)ds.
\end{align*}
For any $0\leq t_0\leq T,$ we get
\begin{align*}
  \sup_{0\leq t\leq t_0}|\widehat{v}(t,\xi)|
 \leq |\widehat{v_0}|(\xi)+\int_0^{t_0} |\widehat{f}|(s,\xi)ds.
\end{align*}
Taking the $L^1$ norm with respect to $|\xi|^{-1}d\xi$ allows us to conclude that
\begin{align*}
  \int_{\mathbb{R}^3}\sup_{0\leq t\leq t_0}|\widehat{v}(t,\xi)|\xi|^{-1}d\xi
 &\leq \int_{\mathbb{R}^3}|\widehat{v_0}|(\xi)|\xi|^{-1}d\xi+\int_{\mathbb{R}^3}\int_0^{t_0} |\widehat{f}|(s,\xi)|\xi|^{-1}dsd\xi\\
 \nonumber&\leq \|v_0\|_{\chi^{-1}}+\|f\|_{L^1([0,t_0];\chi^{-1})}.
\end{align*}
The first result (3.6) is thus proved.\\
Similarly, taking the $L^1$ norm with respect to $|\xi|^{1}d\xi dt,$ one have,
\begin{align*}
 &\int_0^{t_0}\int_{\mathbb{R}^3}|\widehat{v}|(s,\xi)|\xi| d\xi dt\\
\nonumber\leq &\int_0^{t_0}\int_{\mathbb{R}^3}e^{-\nu t|\xi|^2}|\widehat{v_0}|(\xi)|\xi|d\xi dt+
 \int_0^{t_0}\int_{\mathbb{R}^3}\int_0^te^{-\nu (t-s)|\xi|^2}|\widehat{f}|(s,\xi)|\xi|dsd\xi dt\\
 \nonumber\leq & \int_{\mathbb{R}^3}\left(\int_0^{t_0}e^{-\nu t|\xi|^2}|\xi|^{2}dt\right)|\widehat{v_0}(\xi)|\xi|^{-1}d\xi\\
 \nonumber&+
\int_0^{t_0}\int_{\mathbb{R}^3}\left(\int_s^{t_0}e^{-\nu (t-s)|\xi|^2}|\xi|^{2}dt\right) |\widehat{f}|(s,\xi)|\xi|^{-1}dsd\xi \\
\nonumber\leq & \frac{1}{\nu}(\int_{\mathbb{R}^3}|\widehat{v_0}(\xi)|\xi|^{-1}d\xi+\|f\|_{L^1([0,t_0];\chi^{-1})}).
\end{align*}
This implies that
\begin{align}\nu \|v\|_{L^1([0,t_0];\chi^{1})}
  \leq \|v_0\|_{\chi^{-1}}+\|f\|_{L^1([0,t_0];\chi^{-1})}.
\end{align}
The inequalities (3.6) and (3.8) lead to (3.7).\\

Now consider two solutions $u_1$ and $u_2$ with the same initial data $u_0,$ and
assume that
\begin{align*}
u_i \in C([0,T];\chi^{-1})\cap L^1([0,T];\chi^1)\cap L^{\infty}([0,T];\chi^{-1}),~~i=1,2.
\end{align*}
Let $w=u_1-u_2,$ we note that $w$ satisfies
 \begin{align}
\left\{
\begin{array}{l}
\partial_tw-\nu \triangle w=Q(w,u_1)+Q(u_2,w),  \\[1ex]
w|_{t=0}=0.
\end{array}
\right.
\end{align}
Thanks to the inequality (3.4), we have that $$Q(w,u_1)+Q(u_2,w)\in L^1([0,T];\chi^{-1}),$$ and that
\begin{align}
  &\|Q(w,u_1)+Q(u_2,w)\|_{L^1([0,t_0];\chi^{-1})}\\
\nonumber\leq &C\|w\|_{L^{\infty}([0,t_0];\chi^{-1})}^{\frac{1}{2}}
(\|u_1\|_{L^{\infty}([0,t_0];\chi^{-1})}+\|u_2\|_{L^{\infty}([0,t_0];\chi^{-1})})^{\frac{1}{2}}\\
\nonumber&\times\|w\|_{L^1([0,t_0];\chi^{1})}^{\frac{1}{2}}(\|u_2\|_{L^1([0,t_0];\chi^{1})}+\|u_2\|_{L^1([0,t_0];\chi^{1})})^{\frac{1}{2}}.
\end{align}
By virtue of Lemma \ref{zhongyao}, we infer that for any $0
\leq t_0 \leq T,$
\begin{align}
  &\|w\|_{L^\infty([0,t_0];(\chi^{-1})}+\nu \|w\|_{L^1([0,t_0];\chi^{1})}\\
\nonumber  \leq &2\|Q(w,u_1)+Q(u_2,w)\|_{L^1([0,t_0];\chi^{-1})}\\
\nonumber \leq& C\|w\|_{L^{\infty}([0,t_0];\chi^{-1})}^{\frac{1}{2}}
(\|u_1\|_{L^{\infty}([0,t_0];\chi^{-1})}+\|u_2\|_{L^{\infty}([0,t_0];\chi^{-1})})^{\frac{1}{2}}\\
\nonumber &\times\|w\|_{L^1([0,t_0];\chi^{1})}^{\frac{1}{2}}
(\|u_2\|_{L^1([0,t_0];\chi^{1})}+\|u_2\|_{L^1([0,t_0];\chi^{1})})^{\frac{1}{2}}\\
\nonumber\leq & \varepsilon(\|u_1\|_{L^{\infty}([0,t_0];\chi^{-1})}+\|u_2\|_{L^{\infty}([0,t_0];\chi^{-1})})\|w\|_{L^{\infty}([0,t_0];\chi^{-1})}\\
\nonumber
&+C_{\varepsilon}(\|u_1\|_{L^1([0,t_0];\chi^{1})}+\|u_2\|_{L^1([0,t_0];\chi^{1})})\|w\|_{L^1([0,t_0];\chi^{1})}.
\end{align}
Choosing $\varepsilon>0$ such that
\begin{align*}
  \varepsilon(\|u_1\|_{L^{\infty}([0,T];\chi^{-1})}+\|u_2\|_{L^{\infty}([0,T];\chi^{-1})})\leq \frac{1}{2},
\end{align*}
then there exists a positive number $\delta$  satisfying $0<\delta\leq T $ and
\begin{align*}
  C_{\varepsilon}(\|u_1\|_{L^1([0,\delta];\chi^{1})}+\|u_2\|_{L^1([0,\delta];\chi^{1})})\leq \frac{\nu}{2}.
\end{align*}
We then infer that
\begin{align*}
&\|w\|_{L^\infty([0,\delta];\chi^{-1})}+\nu \|w\|_{L^1([0,\delta];\chi^{1})}\\\leq
&\varepsilon(\|u_1\|_{L^{\infty}([0,\delta];\chi^{-1})}+\|u_2\|_{L^{\infty}([0,\delta];\chi^{-1})})\|w\|_{L^{\infty}([0,\delta];\chi^{-1})}\\
\nonumber
&+C_{\varepsilon}(\|u_1\|_{L^1([0,\delta];\chi^{1})}+\|u_2\|_{L^1([0,\delta];\chi^{1})})\|w\|_{L^1([0,\delta];\chi^{1})}\\
\nonumber\leq &\frac{1}{2}(\|w\|_{L^\infty([0,\delta];\chi^{-1})}+\nu \|w\|_{L^1([0,\delta];\chi^{1})}).
\end{align*}
This implies that $w(t)=0$, $0\leq t\leq \delta.$ Basic connective argument then yields uniqueness on $[0,T].$\\
\vspace*{1em}

Theorem \ref{main theo} is thus proved up to the blow-up criterion. Assume that we have a solution of the system $(GNS_\nu)$ on a time interval $[0,T]\,(T<\infty)$ such that
\begin{align}
\int_0^T\|u\|_{\chi^0}^2dt<\infty .
\end{align}
We claim that the lifespan $T_{u_0}$ of $u$ is greater than $T.$  Indeed, thanks to Lemma \ref{zhongyao}, the inequalities (3.4) and (3.12), we have
\begin{align*}
\int_{\mathbb{R}^3} \sup_{0\leq t\leq T}|\widehat{v}|(t,\xi)|\xi|^{-1}d\xi
&\leq \|u_0\|_{\chi^{-1}}+\|Q(u,u)\|_{L^1([0,T];\chi^{-1})}\\
\nonumber&\leq \|u_0\|_{\chi^{-1}}+C\|u\|_{L^2[0,T];\chi^{0})}^2\\
\nonumber& < \infty.
\end{align*}
Thus, a positive number $\rho$ exists such that
\begin{align*}
\forall t\in [0,T],\int_{|\xi|\geq \rho}|\xi|^{-1}|\widehat{u}(t,\xi)|(\xi)d\xi\leq \frac{\nu}{2^{\frac{5}{2}}C_0}.
\end{align*}

The condition (3.3) now implies that for any $t\in[0,T],$ the lifespan for a solution of $(GNS_\nu)$ with initial data $u(t)$ is bounded from below by a positive real number C which is independent of $t.$ Thus $T_{u_0}>T,$ and the whole of Theorem \ref{main theo} is now proved.\\
\\

As the system $(NS_\nu)$ is a particular case of the system $(GNS_\nu),$ we only need to  show that the Navier-Stokes system is well-posed globally in time for $\|u_0\|_{\chi^{-1}}<\nu.$ This result has been obtained in \cite{Zhen Lei} by mollifying initial date. Here, however, due to Theorem \ref{main theo}, for any $u_0\in \chi^{-1}$ there already exists a unique local solution $u(x,t)$ on some internal $[0,T^*),$ thus we can directly obtain estimates of the solution $u$ instead of approximate solutions \cite{Zhen Lei}.\\

Taking the Fourier transform of the system $(NS_\nu)$, one has
\begin{align}\label{*}
&\partial_t\widehat{u}(t,\xi)+\nu |\xi|^2\widehat{u}(t,\xi)\\
\nonumber=&i\int_{\mathbb{R}^3}\widehat{u}(\eta)\otimes\widehat{u}(\xi-\eta)d\eta \cdot \xi -i\left(\frac{1}{|\xi|^2}\xi^{T}\cdot \int_{\mathbb{R}^3}\widehat{u}(\eta)\otimes\widehat{u}(\xi-\eta)d\eta \cdot \xi\right)\xi,
\end{align}
The condition $div\,u=0$ implies that $ \xi \cdot \widehat{u}=0.$ \\

Letting $(\ref{*})\cdot \overline{\widehat{u}}+
\widehat{u}\cdot\overline{(\ref{*})},$ we get
 \begin{align*}
   &\partial_t|\widehat{u}|^2(t,\xi)+2\nu|\xi|^2|\widehat{u}|^2(t,\xi)\\
   =&i\left[\xi^{T}\cdot\int_{\mathbb{R}^3}\widehat{u}(\eta)\otimes\widehat{u}(\xi-\eta)d\eta \cdot \bar{\hat{u}}-\hat{u}^{T}\cdot\int_{\mathbb{R}^3}\bar{\widehat{u}}(\eta)\otimes\bar{\widehat{u}}
   (\xi-\eta)d\eta \cdot \xi\right].
 \end{align*}
For any positive $\varepsilon,$ we have
\begin{align*}
\partial_t(|\widehat{u}|^2+\varepsilon)^{\frac{1}{2}}
=\frac{\partial_t|\widehat{u}|^2}{2(|\widehat{u}|^2+\varepsilon)^{\frac{1}{2}}}.
\end{align*}
Thus, integrating with respect to $t$ gives
\begin{align*}
  &(|\widehat{u}|^2(t,\xi)+\varepsilon)^{\frac{1}{2}}+\nu\int_0^t\frac{|\xi|^2|\widehat{u}|^2(s,\xi)}
  {(|\widehat{u}|^2(s,\xi)+\varepsilon)^{\frac{1}{2}}}ds\\
  \leq &(|\widehat{u_0}|^2(\xi)+\varepsilon)^{\frac{1}{2}}+
  \int_0^t\int_{\mathbb{R}^3}\frac{|\widehat{u}|(s,\xi)|\xi||\widehat{u}|(s,\eta)|\widehat{v}|(s,\xi-\eta)}
  {(|\widehat{u}|^2(s,\xi)+\varepsilon)^{\frac{1}{2}}}d\eta ds.
\end{align*}
Letting $\varepsilon$ tend to zero gives
\begin{align*}
  |\widehat{u}|(t,\xi)+\nu\int_0^t  |\xi|^2\ |\widehat{u}|(s,\xi)ds
  \leq |\widehat{u_0}(\xi)|+\int_0^t \int_{\mathbb{R}^3}|\xi||\widehat{u}|(s,\eta)|\widehat{u}|(s,\xi-\eta)d\eta ds.
\end{align*}
Taking the $L^1$ norm with respect to $|\xi|^{-1}d\xi,$ we have
\begin{align*}
 &\int_{\mathbb{R}^3}|\widehat{u}|(t,\xi)|\xi|^{-1}d\xi+\nu \int_0^t\int_{\mathbb{R}^3}|\widehat{u}|(s,\xi)|\xi|d\xi ds\\
 \leq&\int_{\mathbb{R}^3}|\widehat{u_0}(\xi)|\xi|^{-1}d\xi+\|u\|_{L^2([0,t];\chi^{0})}^2\\
 \leq&\|u_0\|_{\chi^{-1}}+
 \|u\|_{L^\infty([0,t];\chi^{-1})}\|u\|_{L^1([0,t];\chi^{1})}.
\end{align*}
If $\|u_0\|_{\chi^{-1}}< \nu,$ then $\|u(t)\|_{\chi^{-1}}< \nu$ at least for a very short time interval $[0,\delta].$ Consequently, on such a time interval, we have
$$\|u(t)\|_{\chi^{-1}}\leq \|u(0)\|_{\chi^{-1}}< \nu.$$
The basic continuity argument yields that
\begin{align*}
\|u(t)\|_{\chi^{-1}}\leq \|u(0)\|_{\chi^{-1}}< \nu,
\end{align*}
for all $t\in[0,T^*)$ (see the inequalities (2.5)-(2.6) in \cite{Zhen Lei}).\\
We then derive that\begin{align*}
\|u(t)\|_{\chi^{-1}}+(\nu-\|u(0)\|_{\chi^{-1}})\int_0^t\|u(s)\|_{\chi^1} ds\leq\|u_0\|_{\chi^{-1}},
\end{align*}
for all $t\in[0,T^*)$.\\
Thanks to Proposition \ref{chazhi}, one has
\begin{align*}
  \int^{T^*}_0 \|u(t)\|_{\chi^0}^2dt\leq\frac{\|u_0\|_{\chi^{-1}}^2}{\nu-\|u(0)\|_{\chi^{-1}}}<\infty.
\end{align*}
According to Theorem \ref{main theo} again, this implies that $T^*=\infty,$ and the whole Theorem \ref{main coro} is proved.\\
\vspace*{1em}

\noindent\textbf{Acknowledgements}. This work was partially
supported by NNSFC (No. 11271382), RFDP (No. 20120171110014), and
the key project of Sun Yat-sen University (No. c1185).

\phantomsection

\end{document}